\documentclass[11pt]{amsart}
\usepackage{amssymb,amsthm,amsfonts,amsmath}

\newcommand{\Q}{{\mathbb Q}}

\newcommand{\Z}{{\mathbb Z}}

\newcommand{\Proj}{{\mathbf P}}

\newcommand{\EE}{{\mathcal E}}
\newcommand{\C}{{\mathbb C}}
\newcommand{\card}{{{\operatorname{card}}}}

\newcommand{\F}{{\mathbb F}}

\newcommand{\W}{{\mathcal{W}}}

\newcommand{\comment}[1]{}
\newtheorem{theorem}{Theorem}[section]
\newtheorem{lemma}[theorem]{Lemma}

\newtheorem{prop}[theorem]{Proposition}

\theoremstyle{definition}
 
\newtheorem{definition}{Definition}
\theoremstyle{remark}
\newtheorem{rem}{Remark$\!\!$}        
\begin{document}
\title[Function Fields of Varieties over $\C$] {Hilbert's
  Tenth Problem for Function Fields of Varieties over $\C$ }

\author{Kirsten Eisentr\"ager} 

\address{School of Mathematics, Institute for Advanced Study,
  Einstein Drive, Princeton, NJ 08540, USA}

\email{eisentra@ias.edu}
\begin{abstract}
  Let $K$ be the function field of a variety of dimension $\geq 2$ over
  an algebraically closed field of characteristic zero. Then Hilbert's
  Tenth Problem for $K$ is undecidable.  This generalizes the result
  by Kim and Roush from 1992 that Hilbert's Tenth Problem for the
  purely transcendental function field ${\mathbb{C}}(t_1,t_2)$ is
  undecidable.
\end{abstract}
\thanks{Partially supported by the
  National Science Foundation under agreement No.\ DMS-0111298. Any
  opinions, findings and conclusions expressed in this article are
  those of the author and do not necessarily reflect the views of the
  National Science Foundation.}

\maketitle
\section{Introduction}
Hilbert's Tenth Problem in its original form was to find an algorithm
to decide, given a polynomial equation $f(x_1,\dots,x_n)=0$ with
coefficients in the ring $\Z$ of integers, whether it has a solution
with $x_1,\dots,x_n \in \Z$.  Matijasevi{\v{c}}~\cite{Mat70} proved
that no such algorithm exists, {\it i.e.}\ that Hilbert's Tenth
Problem is undecidable. Since then various analogues of this problem
have been studied by considering polynomial equations with
coefficients and solutions over some other commutative ring $R$.
Perhaps the most important unsolved question in this area is Hilbert's
Tenth Problem over the field of rational numbers.  However, there has
been recent progress by Poonen~\cite{Poo03} who proved undecidability
for large subrings of $\Q$. The function field analogue, namely
Hilbert's Tenth Problem for the function field $k$ of a curve over a
finite field, is undecidable.  This was proved by Pheidas for
$k=\F_q(t)$ with $q$ odd, and by Videla \cite{Vi94} for $\F_q(t)$ with
$q$ even.  Shlapentokh~\cite{Sh96,Sh2000} generalized Pheidas' result
to finite extensions of $\F_q(t)$ with $q$ odd and to certain function
fields over possibly infinite constant fields of odd characteristic,
and the remaining cases in characteristic $2$ are treated in
\cite{Eis2002}.  Hilbert's Tenth Problem is also known to be
undecidable for several rational function fields of characteristic
zero: In 1978 Denef proved the undecidability of Hilbert's Tenth
Problem for rational function fields $K(T)$ over formally real fields
$K$ \cite{Den78}, and he was the first to use rank one elliptic curves
to prove undecidability. Kim and Roush \cite{KR92} showed that the
problem is undecidable for the purely transcendental function field
$\C(t_1,t_2)$ and their approach also used rank one elliptic curves.
In this paper we will generalize the result by Kim and Roush. We prove
that Hilbert's Tenth Problem for function fields of surfaces over $\C$
is undecidable.

More generally, in the first part of this paper we prove the following:
\begin{theorem}\label{Theorem}
  Let $K$ be the function field of a surface over an algebraically
  closed field $k$ of characteristic zero. Then Hilbert's Tenth
  Problem for $K$ is undecidable.
\end{theorem}

In Hilbert's Tenth Problem we need to input the coefficients of an
equation into a Turing machine, so in our situation we need to
restrict the coefficients of the equations. We will show: There exist
elements $z_1,z_2 \in K$ which generate an extension of transcendence
degree $2$ over $k$ such that Hilbert's Tenth Problem for $K$ with
coefficients in $\Z[z_1,z_2]$ is undecidable. For simplicity we will
just refer to this as Hilbert's Tenth Problem for $K$.

Theorem~\ref{Theorem} can also be stated more geometrically: Let $X$
be a smooth projective surface over $k$. There is no algorithm that
determines whether dominant morphisms from varieties to $X$ admit
rational sections.

We will prove Theorem~\ref{Theorem} by constructing a diophantine
model of the integers with addition and multiplication over $K$. To do
this we will construct a diophantine model of $\Z \times \Z$ with
certain relations and then show that we can define the integers
with addition and multiplication inside this model. In the first part
of this paper we will for simplicity assume that $K$ is the function
field of a surface over $\C$.  The proof is exactly the same when we
replace $\C$ by an arbitrary algebraically closed field of
characteristic zero. We use the fact that our field has characteristic
zero in the construction of rank one elliptic curves in
Section~\ref{elliptic}, and the fact that we are working with function
fields over an algebraically closed field is used in Section~\ref{div}
where we apply the Tsen-Lang Theorem.


The same approach as in Theorem~\ref{Theorem} also generalizes to
higher transcendence degree. In Section~\ref{generalization} we prove:
\begin{theorem}\label{general}
  Let $K$ be the function field of a variety of dimension $\geq2$ over
  an algebraically closed field $k$ of characteristic zero. Then
  Hilbert's Tenth Problem for $K$ is undecidable.
\end{theorem}
Again, we will show: There exist elements $z_1,z_2 \in K$
which generate an extension of transcendence degree $2$ over $k$ such
that Hilbert's Tenth Problem for $K$ with coefficients in
$\Z[z_1,z_2]$ is undecidable.
\section{The model}\label{TheModel}
First we will define two notions that will appear frequently in
the remainder of this paper.
\begin{definition}
  1. If $R$ is a commutative ring, a {\em diophantine equation over
    $R$} is an equation $P(x_1,\dots,x_n)=0$ where $P$ is a polynomial
  in the variables $x_1, \dots, x_n$ with coefficients in
  $R$.\\
  2. A subset $S$ of $R^k$ is {\em diophantine over R} if there exists a
  polynomial $P(x_1,\dots,x_k,$ $y_1,\dots,y_m) \in
  R[x_1,\dots,x_k,y_1,\dots,y_m]$ such that \[S=\{(x_1, \dots, x_k)
  \in R^k: \exists y_1, \dots,y_m \in R,\;
  (P(x_1,\dots,x_k,y_1,\dots,y_m)=0)\}.\]
\end{definition}

Let $K$ be a finite extension of $\C(t_1,t_2)$.  We will define a diophantine
model of the structure $\mathcal{S}= \langle \,\Z \times
\Z,+,\mid\,,\mathcal{Z},\mathcal{W}\,\rangle$ in $K$. Here
$+$ denotes the usual component-wise addition of pairs of integers,
$\mid$ re\-presents a relation which satisfies
\[(n,1)\mid(m,s)\Leftrightarrow m = ns,\]
and $\mathcal{Z}$ is a unary predicate which is interpreted as
\[
\mathcal{Z}(n,m)\Leftrightarrow m=0.
\]
The predicate $\W$ is interpreted as 
\[
\W((m,n),(r,s)) \Leftrightarrow m=s \land n=r.
\]
A {\em diophantine model} of $\mathcal{S}$ over $K$ is a diophantine
subset $S \subseteq K^n$ equipped with a bijection $\phi: \Z \times \Z
\to S$ such that under $\phi$, the graphs of addition, $\mid$,
$\mathcal{Z}$, and $\mathcal{W}$ in $\Z \times \Z$ 
 correspond to diophantine subsets of $S^3, S^2, S,$ and $S^2$, respectively.

In this section we will show that constructing such a model is
sufficient to prove undecidability of Hilbert's Tenth Problem for $K$.

First we can show the following
\begin{prop}
The relation $\W$ can be defined entirely in terms of the
other relations.
\end{prop}
\begin{proof}
It is enough to verify that
\[
\W((a,b),(x,y)) \Leftrightarrow (1,1)\mid((x,y)+(a,b)) \land
        (-1,1)\mid((x,y)-(a,b)).\]
\end{proof}

As Pheidas and Zahidi \cite{PheZa00} point out we can existentially
define the integers with addition and multiplication inside
$\mathcal{S}= \langle \Z \times
\Z,+,\mid\,,\mathcal{Z},\W\rangle$, so $\mathcal{S}$ has an
undecidable existential theory:

\begin{prop}
The structure $\mathcal{S}$ has an undecidable existential theory.
\end{prop}
\begin{proof}
  We interpret the integer $n$ as the pair $(n,0)$. The set $\{(n,0):n
  \in \Z\}$ is existentially definable in $\mathcal{S}$ through the
  relation $\mathcal{Z}$. Addition of integers $n,m$ corresponds to the
  addition of the pairs $(n,0)$ and $(m,0)$. To define multiplication
  of the integers $m$ and $r$, note that $n=mr$ if and only if
  $(m,1)\mid(n,r)$, hence $n=mr$ if and only if
\[
\exists\, a,b \,\,:((m,0)+(0,1))\mid((n,0)+(a,b))\land
\W((a,b),(r,0)).
\]
Since the positive existential theory of the integers with addition
and multiplication is undecidable, $\mathcal{S}$ has an undecidable
existential theory as well.
\end{proof}
The above proposition shows that in order to prove
Theorem~\ref{Theorem} it is enough to construct a diophantine model of
$\mathcal{S}$ over $K$.  In Sections~\ref{definable} and ~\ref{div} we
will construct this model.

\section{An existentially definable set isomorphic to $\Z \times \Z$}
\label{definable}

\subsection{Generating elliptic curves of rank one}\label{elliptic}
Let $K$ be a finite extension of $\C(t_1,t_2)$.

Our first task is to find a diophantine set $A$
over $K$ which is isomorphic to $\Z \times \Z$ as a set.  In their
undecidability proof for $\C(t_1,t_2)$ \cite{KR92} Kim and Roush
obtain such a set by using the $\C(t_1,t_2)$-rational points on two
elliptic curves which have rank one over $\C(t_1,t_2)$. We need to
construct two elliptic curves which have rank one over $K$.

These can be obtained from a theorem by Moret-Bailly~\cite{MB01}. In
his theorem he uses the following notation: Let $k$ be a field of
characteristic zero. Let $C$ be a smooth projective geometrically
connected curve over $k$ with function field $F$. Let $Q$ be a finite
nonempty set of closed points of $C$.  Let $E:y^2 = x^3 + ax +b$ be an
elliptic curve over $k$ with $b \neq 0$.  In \cite{MB01} Moret-Bailly
also introduces another curve $\Gamma$, but for our application of his
theorem we only have to consider the special case where $\Gamma=E$.
\begin{definition}
Let $k,C,F,E,Q$ be as above.
Let $g:C\to \Proj^1_k$ be a non-constant $k$-morphism
corresponding to an injection $k(T) \hookrightarrow F$ sending $T$ to $g$.
We say that $g$ is {\em admissible} for $E$ (and $Q$)  if
\begin{enumerate}
\item $g$ has only simple branch points.
\item $g$ is \'etale above 0 and the branch points of $E$.
\item Every point of $Q$ is a pole of $g$.
\end{enumerate}
\end{definition}
\begin{rem}
  In \cite{MB01} Moret-Bailly asserts that, given $C,E,Q$, we can
  always find an admissible morphism $g$. If $g$ is admissible for
  $E$, then for all but finitely many $\lambda \in k^{*}$, $\lambda g$
  is still admissible.
\end{rem}
Now we can state Moret-Bailly's theorem.
\begin{theorem}\cite[Theorem 1.5]{MB01}\label{MB}
  Let $k,C,F,E, Q$ be as above.  Let $f \in F$ be admissible for
  $E,Q$.  Let $\EE_{\lambda f}$ be the elliptic curve $$\EE_{\lambda
    f}: ((\lambda f)^3 + a(\lambda f) +b)\,y^2 = x^3 + ax + b.$$ Then
  the natural homomorphism $\EE(k(T))\hookrightarrow \EE_{\lambda f}(F)$
  induced by the inclusion $k(T) \hookrightarrow F$ that sends $T$ to
  $\lambda f$ is an isomorphism for infinitely many $\lambda \in \Z$.
\end{theorem}
If we assume in addition that $E$ has no complex multiplication, then by
a theorem of Denef \cite{Den78}, $\EE(k(T))$ has rank one with
generator $(T,1)$ (modulo 2-torsion). So Theorem~\ref{MB} gives us
infinitely many elliptic curves that have rank one over $F$: Given
$f$ and $\lambda$ as in the theorem we obtain an elliptic curve
$\EE_{\lambda f}$ which has rank one over $F$ with generator $(\lambda
f,1)$ (modulo 2-torsion).

We can use Moret-Bailly's theorem to prove the following theorem:
\begin{theorem}\label{rank2curves}
  Let $K$ be a finite extension of $\C(t_1,t_2)$.  Let $E: y^2 = x^3 +
  ax + b$ be an elliptic curve with $a,b \in \C$, $b\neq 0$, and no
  complex multiplication.  There exist $z_1,z_2 \in K$ such that
  $\C(z_1,z_2)$ has transcendence degree $2$ over $\C$ and such that
  the two elliptic curves $\EE_1: (z_1^3 + az_1+b)\,y^2 = x^3 + ax +
  b$ and $\EE_2: (z_2^3 + az_2+b)\,y^2 = x^3 + ax + b$ have rank one
  over $K$ with generators $(z_1,1)$ and $(z_2,1)$, respectively
  (modulo 2-torsion).
\end{theorem}
\begin{proof}
  Let $k$ be the algebraic closure of $\C(t_2)$ inside $K$. There
  exists a smooth, projective, geometrically connected curve $C$ over
  $k$ whose function field is $K$. Let $Q$ be a finite nonempty set of
  closed points of $C$, and choose an $f \in K$ that is admissible for
  $E$ and $Q$.  Since $f$ is non-constant, $f$ is transcendental over
  $\C(t_2)$.  Now we can apply Theorem~\ref{MB} with $F=K$ and
  $k,C,E,Q, f$ as defined above. By Theorem~\ref{MB} there exists a
  nonzero $\lambda \in \Z$ such that
  $$\EE_{\lambda f}: ((\lambda f)^3 + a (\lambda f) + b )y^2 =x^3 + ax
  +b$$
  has rank one with generator $(\lambda f, 1)$ (modulo
  2-torsion).  By the remark, some integer multiple of $t_2 \cdot
  \lambda f$ will still be admissible for $E$. Let $g$ be such an
  integer multiple. By Theorem~\ref{MB} applied to $g$ there exists a
  nonzero integer $\mu$ such that $$\EE_{\mu g}: ((\mu g)^3 + a (\mu
  g) + b )y^2 =x^3 + ax +b$$
  has rank one with generator $(\mu g, 1)$
  (modulo 2-torsion).  Let $z_1:=\lambda f$, and let $z_2:= \mu g$. To
  complete the proof it remains to show that $\C(z_1,z_2)$ has
  transcendence degree 2 over $\C$. Since $z_2= \nu z_1 t_2$ for some
  nonzero integer $\nu$, it follows that $t_2 \in \C(z_1,z_2)$. As
  pointed out above, the element $f$ is transcendental over $\C(t_2)$,
  and hence the same is true for $\lambda f = z_1$. This shows that
  the transcendence degree of $\C(z_1,z_2)$ over $\C$ is at least 2,
  and since $\C(z_1,z_2) \subseteq K$, which is algebraic over
  $\C(t_1,t_2)$, the transcendence degree must equal 2.
\end{proof}
\subsection{Diophantine definition of $A$}\label{DiophDef}
In the following let $E: y^2 = x^3 + ax + b$ be an elliptic curve with
$a,b \in \C$, $b \neq 0$, no complex multiplication and such that the point
$(0,\sqrt b)$ has infinite order.  Such a curve exists: the curve
$E:y^2 = x^3 + x+ 1$ (496A1 in \cite{Cremona}) has the required
properties. The condition that $(0,\sqrt b)$ has infinite order will
be needed in the proof of Theorem~\ref{almostdiv}. Let $\EE_1$ and
$\EE_2$ be as above. To be able to define a suitable set $A$ which is
isomorphic to $\Z \times \Z$ we need to work in an algebraic extension
$L$ of $K$. Let $L:=K(h_1,h_2)$, where $h_i$ is defined by
$h_i^2=z_i^3 + az_i+b$, for $i=1,2$. To prove undecidability for $K$
it is enough to prove that the existential theory of $L$ in the
language $\langle\, L,+,\cdot \; ; 0,1,z_1,z_2,h_1,h_2,S\,\rangle$ is
undecidable, where $S$ is a predicate for the elements of the subfield
$K$ \cite[Lemma 1.9]{PheZa00}. So from now on we will work with
equations over $L$.

Over $L$ both $\EE_1$ and $\EE_2$ are isomorphic to $E$. There is an
isomorphism between $\EE_1$ and $E$ that sends $(x,y) \in \EE_1$ to
the point $(x,h_1y)$ on $E$. Under this isomorphism the point
$(z_1,1)$ on $\EE_1$ corresponds to the point $P_1:=(z_1,h_1)$ on $E$.
Similarly there is an isomorphism between $\EE_2$ and $E$ that sends
the point $(z_2,1)$ on $\EE_2$ to the point $P_2:=(z_2,h_2)$ on $E$.

The elliptic curve $E$ is a projective variety, but any projective
algebraic set can be partitioned into finitely many affine algebraic
sets, which can then be embedded into a single affine algebraic
set. This implies that the set $E(L)$ is diophantine over $L$, since
we can take care of the point at infinity $\mathbf{O}$ of $E$.

We will show that the set of points $\Z P_1 \times \Z P_2 \subseteq
E(L)$ is existentially definable in our language, because we have a
predicate for the elements of $K$: Since $\EE_1$ has $2$-torsion, we
first give a diophantine definition of $2 \cdot \, \Z P_1$:
\[
P \in 2 \cdot \,\Z P_1 \Leftrightarrow \exists \,x,y \in K \,\, (z_1^3 +
a z_1 + b)\,y^2 = x^3 + ax + b \land P = 2 \cdot (x, h_1y)
\]
This set is diophantine, because we can express that $P = 2\cdot Q$ by
diophantine equations.  Then $\Z P_1$ can be defined as
\[
P \in \Z P_1 \Leftrightarrow (P \in 2\cdot \, \Z  P_1) {\mbox{ or }}
(\exists \,Q \in 2 \cdot \, \Z P_1\mbox{ and } P = Q + P_1)
\]
Similarly we have a diophantine definition for $\Z P_2$. Hence the
cartesian product $\Z P_1 \times \Z P_2 \subseteq E(L)$ is
existentially definable, since addition on $E$ is existentially
definable.
\subsection{Existential definition of $+$ and $\mathcal{Z}$}\label{relations}
The unary relation $\mathcal{Z}$ is existentially definable, since
this is the same as showing that the set $\Z P_1$ is diophantine,
which was done above.  Addition of pairs of integers corresponds to
addition on the cartesian product of the elliptic curves $\EE_i$ (as
groups), hence it is existentially definable.  Since $\W$ can be
defined in terms of the other relations, it remains to define the
divisibility relation $\mid$\,. This is done in the next section.

\section{Existential definition of
  $(m,1)\mid(n,r)$ }
\label{div}
Now we will show how to existentially define the relation $\mid$ among
pairs of integers in $L$.
In the following $x(P)$ will denote the $x$-coordinate of a
point $P$ on $E$, and $y(P)$ will denote the $y$-coordinate of $P$.
Let $\alpha := [L:\C(z_1,z_2,h_1,h_2)]$. To give an existential
definition of $\mid$ we will use the fact that $(m,1)\mid(n,r)
\Leftrightarrow (m,1)\mid(kn,kr)$ for $k \neq 0$.
\begin{theorem}\label{almostdiv}
  There exists a finite set $U \subseteq \Z$ such that for all $m \in
  \Z-U$ we have: for all $n,r \in \Z$
\begin{align*} (m,1)\mid(n,r) \Leftrightarrow\\ 
  (\,\exists \, y_0,z_0 \in L^{*}\, \,\, &x(nP_1+rP_2)\, y_0^2 +
    x(mP_1+P_2)\,z_0^2 =1 \\ \land \,\, \exists\, y_1,z_1 \in
  L^{*} \,\,\,
  &x(2nP_1+2rP_2)\, y_1^2 + x(mP_1+P_2)\,z_1^2 =1 \\
  \cdots\\ \land \,\, \exists\,
  y_{\alpha},z_{\alpha} \in L^{*} \,\,\,
  &x(2^{\alpha}nP_1+2^{\alpha}rP_2)\, y_{\alpha}^2 +
  x(mP_1+P_2)\,z_{\alpha}^2 =1 \,)\,.
\end{align*}
\end{theorem}
\begin{proof}
  For the first implication, assume that $(m,1)\mid(n,r)$, {\it i.e.}\
  $n=mr$.  Then both $ x(nP_1+rP_2)=x(r(mP_1+P_2))$ and $x(mP_1+P_2)$
  are elements of $\C(x(mP_1+P_2), y(mP_1+P_2))$, which has
  transcendence degree one over $\C$. This means that we can apply the
  Tsen-Lang Theorem (Theorem~\ref{Tsen-Lang} in the appendix) to the
  quadratic form
\[
x(nP_1+rP_2)y^2 + x(mP_1+P_2)z^2 -w^2
\]
to conclude that there exists a nontrivial zero $(y,z,w)$ over
$\C(x(mP_1+P_2), y(mP_1+P_2))$. From the theory of quadratic forms it
follows that there exists a nontrivial zero $(y,z,w)$ with $y\cdot z
\cdot w\neq 0$.  The same can be done for the other equations.

For the other direction, suppose that $n \neq mr$ and assume by
contradiction that all $\alpha+1$ equations are satisfied.  We will
proceed with the proof in four steps.
\\
{\bf Claim 1.} There exists a finite set $U \subseteq \Z$ such that
for all $m \in \Z-U$ there exists a discrete valuation $w_m:L^{*}
\twoheadrightarrow \Z$ such that $w_m(x(mP_1+P_2))=1$ and such that
$w_m(x(knP_1+krP_2))=~0$ for $k=1,2,4,\dots,2^{\alpha}$.
\\
{\bf Proof of Claim 1:} Fix $m \in \Z$. Let $P_2' = mP_1 + P_2
=(z_2',h_2')$. Let $F:= \C(z_1,z_2,h_1,h_2)$. Then
\[F=\C(z_1,z_2,h_1,h_2)=\C(z_1,h_1,z_2',h_2')=
\C(x(P_1),y(P_1),x(P_2'),y(P_2')),\] since $z_2 = x(P_2'-mP_1)$ and
$z_2'=y(P_2'-mP_1)$.  Let $v_m:F^{*}\twoheadrightarrow \Z$ be a
discrete valuation which extends the discrete valuation $v$ of
$\C(z_1,h_1)(z_2')$ associated to $z_2'$. The valuation $v$ is the
discrete valuation that satisfies $v(\gamma)=0$ for all $\gamma \in
\C(z_1,h_1)$ and $v(z_2')=1$. In $F$ the element $z_2'$ is still a
uniformizer.  Suppose $v_m$ does not ramify in $L$. Let $w_m$ be an
extension of $v_m$ to $L$. We can consider $L/F$ as an extension of
algebraic function fields of transcendence degree one by taking the
constant field of $F$ to be $\C(z_1,h_1)$, and $v_m$ is a valuation
which is trivial on $\C(z_1,h_1)$. We can apply Theorem~\ref{ramified}
to conclude that only finitely many $v_m$ ramify in $L$. Let \[U:=
\{m\in \Z : v_m \text{ ramifies in }L\}.\] Then $U$ is a finite set.
Let $m \in \Z-U$, and let $\ell$ be the residue field of $w_m$. Since
the residue field of $v_m$ is just $\C(z_1,h_1)$, it follows that
$\ell$ is a finite extension of $\C(z_1,h_1)$ with $[\ell:\C(z_1,h_1)]
\leq \alpha$.  Let $s:= n-mr$.  By assumption $s\neq 0$.  The
equations, rewritten in terms of $P_2'$ and $s$ become
\begin{align*}
\exists \,\,y_0,z_0 \in L^{*} \,\,\, &x(sP_1+rP_2') \,y_0^2 + x(P_2')\,z_0^2
=1\\ \land\,\, \exists\,\, y_1,z_1 \in L^{*} \,\,\,
&x(2sP_1+2rP_2')\, y_1^2 + x(P_2')\,z_1^2 =1 \\
 \cdots\\ \land \,\, \exists\,\,
y_{\alpha},z_{\alpha} \in L^{*} \,\,\,
&x(2^{\alpha}sP_1+2^{\alpha}rP_2')\, y_{\alpha}^2 +
x(P_2')\,z_{\alpha}^2 =1 .
\end{align*}
By our choice of $w_m$, we have $w_m(x(P_2'))=1$. 
The residue field $\ell$ of $w_m$ is an extension of $\C(z_1, h_1)$
and the image of $x(sP_1+rP_2')$ in the residue field $\ell$ is
$x(s(z_1,h_1)+r(0,\pm \sqrt b))$.  We can show that this
$x$-coordinate cannot be zero, which will imply that
$w_m(x(nP_1+rP_2))=w_m(x(sP_1+rP_2')) =0$: The point $(z_1,h_1) \in
E(\C(z_1,h_1))$ has infinite order, $z_1$ is transcendental over $\C$,
and all points of $E$ whose $x$-coordinate is zero are defined over
$\C$. Since $s\neq 0$, this implies that $x(s(z_1,h_1)+r(0,\pm \sqrt
b)) \neq 0$. 

The same argument shows that $w_m(x(knP_1+krP_2))=w_m(x(k s P_1+k r
P_2'))=0$ for $k=1,2,4,\dots,2^{\alpha}$ for all $m \in \Z-U$.
\\
{\bf Claim 2.} Denote by $x_{s,r}$ the image of
$x(sP_1+rP_2')=x(nP_1+rP_2)$ in the residue field $\ell$.  The
elements $x_{s,r}\;$, $x_{2s,2r}\;, \dots,
\;x_{2^{\alpha}s,2^{\alpha}r}$ are squares in $\ell$.
\\
{\bf Proof of Claim 2:} This follows immediately from
Lemma~\ref{square} in the appendix.
\\
{\bf Claim 3.} The elements $x_{s,r}\;, \;\dots,\;
x_{2^{\alpha}s\;,\;2^{\alpha}r}$ are not squares in $\C(z_1,h_1)$.
\\
{\bf Proof of Claim 3:} Since $x_{s,r}\;,x_{2s,2r}\;, \cdots,
x_{2^{\alpha}s,2^{\alpha} r} \in \C(z_1,h_1)$, which is the function
field of $E$, we can consider them as functions $E \to
\Proj^1_{\C}$. Then $x_{s,r}$ corresponds to the function on $E$ which
can be obtained as the composition $P \mapsto sP+ r (0, \sqrt
b)\mapsto x(sP+r(0,\sqrt b))$. The $x$-coordinate map is of degree $2$
and has two distinct zeros, namely $(0, \sqrt b)$ and $(0, -\sqrt b)$.
The map $E \to E$ which maps $P$ to $(sP + r (0,\sqrt b))$ is
unramified since it is the multiplication-by-$s$ map followed by a
translation.  Hence the composition of these two maps has $2s^2$
simple zeros. The same argument works for the other functions
$x_{ks,kr}$. So each of the functions $x_{ks,kr}$, for
$k=1,2,4,\dots,2^{\alpha}$ has only simple zeros.  In particular, none
of these functions is a square in $\C(z_1,h_1)$.
\\
{\bf Claim 4.} The images of $x_{s,r}, \dots,
x_{2^{\alpha}s,2^{\alpha}r}$ in \[V:=[(\ell^*)^2 \cap\,
\C(z_1,h_1)^*]/(\C(z_1,h_1)^*)^2 \text{ are distinct.}\] {\bf Proof of
  Claim 4:} By Claim 2 all elements $x_{ks,kr}$ are in $(\ell^*)^2$.
If $x_{s,r}(P) = 0$ for some point $P$ on $E$, then
\[sP+r(0,\sqrt b) = (0, \sqrt b) \mbox{ or } sP+r(0,\sqrt b) =
(0,-\sqrt b). \] But that implies that
\[2[sP+r(0,\sqrt b)] \neq \mathbf{O},(0, \sqrt b), (0, - \sqrt
b),\] since we picked our curve $E: y^2 = x^3 + ax + b$ such that the
point $(0,\sqrt b)$ has infinite order. Hence neither $2(0,\sqrt b)$
nor $2(0,- \sqrt b)$ can equal $\mathbf O$, $(0,\sqrt b)$ or $(0, -
\sqrt b)$.  This implies that $P$ is neither a zero nor a pole of
$x_{2s,2r}$ and similarly neither a zero nor a pole of
$x_{2^ks,2^kr}$. The same argument shows that a zero of
$x_{2^is,2^ir}$ is neither a zero nor a pole of $x_{2^js,2^jr}$ for $j
>i$.  This implies that it cannot happen that $x_{2^is,2^ir} = f^2
\cdot x_{2^js,2^jr}$ with $f \in \C(z_1,h_1)$, because if, say, $j>i$
and $P$ is a zero of the left-hand-side, then the left-hand-side has a
zero of order $1$ at $P$ while the right-hand-side has a zero of even
order at $P$. Hence all the elements are different in $V=({\ell^*}^2
\cap \C(z_1,h_1)^*)/(\C(z_1,h_1)^*)^2$. This proves the claim.\medskip

But now we have obtained a contradiction: since $[\ell:\C(z_1,h_1)]
\leq \alpha$, the size of $V$ is bounded by $\alpha$ by
Theorem~\ref{Kummer}, so it cannot contain $\alpha+1$ distinct
elements.
This means that for all $m \in \Z -U$ the solvability of the $\alpha
+1$ equations implies that $n=mr$.
\end{proof}

We have seen above that the relation $\mathcal{W}$ can be defined in
terms of the other relations. It turns out that it is convenient to
give an existential definition of $\mathcal{W}$ now and then use it to
give a short proof that $\mid$ has an existential definition.

\begin{prop}
The relation $\mathcal{W}$ is existentially definable.
\end{prop}
\begin{proof}
Let $m_0 \in \Z-U$. Then
\begin{align*}&\mathcal{W}((m,n),(r,s))\\ &\Leftrightarrow
  (1,1)\mid(m+r,n+s) \land (-1,1)\mid(m-r,n-s)\\ &\Leftrightarrow
  (m_0,1)\mid(m_0(m+r),n+s) \land (m_0,1)\mid(-m_0(m-r),n-s) .\end{align*}
Since $m_0$ is a fixed element of $\Z-U$, and since $\Z P_1$ and $\Z
P_2$ are diophantine, the expression
\[(m_0,1)\mid(m_0(m+r),n+s) \land (m_0,1)\mid(-m_0(m-r),n-s)\] is
diophantine in $(m,n)$ and $(r,s)$ by Theorem~\ref{almostdiv}.
\end{proof}
\begin{theorem}
The relation $(m,1)\mid(n,r)$ is existentially definable in $(m,1)$ and
$(n,r)$.
\end{theorem}
\begin{proof}
  Let $m_0$ be as in the above proposition, and let $d$ be a positive
  integer such that $U \subseteq (m_0-d,m_0+d)$.
  Since $n=mr \Leftrightarrow dn+m_0r=dmr+m_0r=(dm+m_0)r$, we have
\[(m,1)\mid(n,r)
\Leftrightarrow (dm+m_0,1)\mid(dn+m_0r,r),\] and we can just work with
that formula instead. So
\[(m,1)\mid(n,r) \Leftrightarrow \exists \,a,b\,
(dm+m_0,1)\mid((dn,r)+m_0(a,b)) \land \mathcal{W}((a,b),(0,r)).\]
This last expression is existentially definable in $(m,1)$ and
$(n,r)$.
Since $m_0 \in \Z -U$, and by our choice of $d$, we have $(dm+m_0) \in
\Z -U$, and we can apply Theorem~\ref{almostdiv}.
\end{proof}
\section{Hilbert's Tenth Problem for function fields of transcendence
  degree $\geq 2$}\label{generalization} In this section we will
generalize Theorem~\ref{Theorem} to function fields of transcendence
degree $\geq 2$ and prove Theorem~\ref{general}. Again, for simplicity
of notation, 
we will prove Theorem~\ref{general} for $k=\C$ . The proof still works
word for word the same when we replace $\C$ by an arbitrary
algebraically closed field of characteristic zero.

Let $K$ be a finite extension of $\C(t_1, t_2,\dots, t_n)$. We want to
prove that Hilbert's Tenth Problem for $K$ is undecidable. In our
proof we will use the following approach. Let $R$ be the algebraic
closure of $\C(t_3,\dots,t_n)$ in $K$. Then $K$ is a finite extension
of $R(t_1,t_2)$. Our argument will follow the proof of the
transcendence degree 2 case with $\C$ replaced by $R$. The only
difference is that we are now working with a transcendence degree 2
extension over a ground field that is no longer algebraically closed.
We only have to modify the parts of the proof of Theorem~\ref{Theorem}
that used the fact that the ground field was algebraically closed.
Theorem~\ref{almostdiv} needed this assumption, because we used the
Tsen-Lang Theorem (Theorem~\ref{Tsen-Lang}) in its proof, but this
theorem is the only part of the proof which relied on the assumption
that we are over an algebraically closed ground field.  As in the
proof of Theorem~\ref{Theorem} we will construct a diophantine model
of the structure $\mathcal{S}= \langle \,\Z \times
\Z,+,\mid\,,\mathcal{Z},\mathcal{W}\,\rangle$ in $K$.  To set up some
notation we will state Theorem~\ref{rank2curves} for higher
transcendence degree.
\begin{theorem}[Theorem~\ref{rank2curves} for higher transcendence degree]
  Let $K$ be a finite extension of $\C(t_1,t_2, \dots, t_n)$, and let
  $R$ be as above.  Let $E: y^2 = x^3 + ax + b$ be an elliptic curve
  with $a,b \in \C$, $b \neq 0$, and no complex multiplication.  There
  exist $z_1,z_2 \in K$ such that $R(z_1,z_2)$ has transcendence
  degree $2$ over $R$ and such that the two elliptic curves $\EE_1:
  (z_1^3 + az_1+b)\,y^2 = x^3 + ax + b$ and $\EE_2: (z_2^3 +
  az_2+b)\,y^2 = x^3 + ax + b$ have rank one over $K$ with generators
  $(z_1,1)$ and $(z_2,1)$, respectively (modulo 2-torsion).
\end{theorem}
\begin{proof}
  The proof goes through word for word as the proof of
  Theorem~\ref{rank2curves} with $\C$ replaced by $R$. 
\end{proof}
As before, we let $L:= K(h_1,h_2)$, where $h_i^2=z_i^3+az_i+b$, and we
will work with the two points $P_1:=(z_1,h_1)$ and $P_2:=(z_2,h_2)$ on
$E$. To prove undecidability for $K$ it is enough to
prove that the existential theory of $L$ in the language $\langle\,
L,+,\cdot \; ; 0,1,z_1,z_2,h_1,h_2,S\,\rangle$ is undecidable, where
$S$ is a predicate for the elements of the subfield $K$ \cite[Lemma
1.9]{PheZa00}. So from now on we will work in this language. We
can use the same argument as in Section~\ref{DiophDef} to show that
the set of points $\Z P_1 \times \Z P_2 \subseteq E(L)$ is
existentially definable in our language. The argument from
Section~\ref{relations} implies that the relations $+$ and
$\mathcal{Z}$ are existentially definable. The only remaining part of
the proof is to show that the divisibility relation $\mid$ is
existentially definable in our language. For this we only have to
prove the analogue of Theorem~\ref{almostdiv}.
\begin{theorem}[Theorem~\ref{almostdiv} for higher transcendence degree]
  There exists a finite set $U \subseteq \Z$ such that for all $m \in
  \Z-U$ we have: for all $n,r \in \Z$
\begin{align*} (m,1)\mid(n,r) \Leftrightarrow\\ 
  (\,\exists \, y_0,z_0 \in L^{*}\, \,\, &x(nP_1+rP_2)\, y_0^2 +
    x(mP_1+P_2)\,z_0^2 =1 \\ \land \,\, \exists\, y_1,z_1 \in
  L^{*} \,\,\,
  &x(2nP_1+2rP_2)\, y_1^2 + x(mP_1+P_2)\,z_1^2 =1 \\
\cdots\\ \land \,\, \exists\,
  y_{\alpha},z_{\alpha} \in L^{*} \,\,\,
  &x(2^{\alpha}nP_1+2^{\alpha}rP_2)\, y_{\alpha}^2 +
  x(mP_1+P_2)\,z_{\alpha}^2 =1 \,)\,.
\end{align*}
\end{theorem}
\begin{proof}
  For the first implication, assume that $(m,1)\mid(n,r)$, {\it i.e.}\ 
  $n=mr$.  Let $K':=\C(x(mP_1+P_2), y(mP_1+P_2))$.  As in the proof of
  Theorem~\ref{almostdiv}, both $ x(nP_1+rP_2)=x(r(mP_1+P_2))$ and
  $x(mP_1+P_2)$ are elements of $K'$, which has transcendence degree
  one over $\C$. This means that we can apply the Tsen-Lang Theorem
  (Theorem~\ref{Tsen-Lang} in the appendix) with $K=K'$ to the quadratic form
\[
x(nP_1+rP_2)y^2 + x(mP_1+P_2)z^2 -w^2
\]
to conclude that there exists a nontrivial zero $(y,z,w)$ over $K'$
and hence that there exists a nontrivial zero $(y,z,w)$ in $K'$ with
$y\cdot z \cdot w\neq 0$. Since $K' \subseteq L$, this produces the
desired solution over $L^{*}$ for the first equation. The
same can be done for the other equations.\\
The other implication in the proof of Theorem~\ref{almostdiv} does not
use the fact that we have a function field over an algebraically
closed field, so we can just work with our finite extension $K$ of
$R(z_1,z_2)$ and with $L$ as defined above, and repeat this part of
the proof word for word with $\C$ replaced by $R$.
\end{proof}
The rest of the proof is word for word the same as in the
transcendence degree 2 case showing that the relation $\mid$ is
existentially definable and thus completing the proof of
Theorem~\ref{general}.
\section{Appendix}\label{auxiliary}
In this section we will state and prove some general theorems which
were used to prove that the relation $\mid$ is
diophantine.
\begin{theorem}\label{ramified}
  Let $L$ and $K$ be function fields of one variable with constant
  fields $C_L$ and $C_K$, respectively, such that $L$ is an extension
  of $K$. If $L$ is separably algebraic over $K$, then there are at
  most a finite number of places of $L$ which are ramified over $K$.
\end{theorem}
\begin{proof}
This theorem is proved on p.~111 of \cite{Deu73} when $C_L \cap K
=C_K$, and the general theorem also follows.
\end{proof}
We also need the following easy lemma.
\begin{lemma}\label{square}
  Let $k$ be a field, and let $v:k^{*} \twoheadrightarrow \Z$ be a
  discrete valuation on $k$. Let $a,b \in k$ with $v(a) =0$ and $v(b)$
  odd.  Suppose that $ax^2 + by^2 =1$ has a solution in $k ^*$. Then
  $a$ is a square in the residue field of $k$.
\end{lemma}
\begin{proof}
  The equation $ax^2 + by^2 =1$ implies that $v(ax^2+ by^2) =0$.  The
  condition $v(a)=0$ implies that $v(ax^2)$ is even. Since $v(b)$ is
  odd, it follows that $v(by^2)$ is odd. Hence $v(ax^2) = 0$ and
  $v(by^2) >0$.  So in the residue field our equation becomes $\bar a
  \cdot \bar x ^2 +0 \equiv 1 \mod v$. This implies that $a$ is a
  square in the residue field.
\end{proof}
\begin{theorem}{Tsen-Lang Theorem.}\label{Tsen-Lang}
Let K be a function field of transcendence degree $j$ over an
algebraically closed field $k$. Let $f_1,\cdots,f_r$ be forms in $n$
variables over $K$, of degrees $d_1, \cdots, d_r$. If
\[
n > \sum_{i=1}^{r}d_i^j
\]
then the system $f_1=\cdots=f_r=0$ has a non-trivial zero in $K^n$.
\end{theorem}
\begin{proof}
  This is proved in Proposition 1.2 and Theorem 1.4 in Chapter 5 of
  \cite{Pfi95}.
\end{proof}
\begin{lemma}\label{Kummer}
  Let $F,G$ be fields of characteristic $\neq 2$, and let $G/F$ be a
  field extension of degree r. Then the cardinality of $V:= [(G^*)^2
  \cap F^*]/(F^*)^2$ is bounded by $r$.
\end{lemma}
\begin{proof}
  The set $V$ is a vector space over $\F_2$. If we have $s$ elements
  of $(G^{*})^2 \cap F^{*}$ whose images in $V$ are linearly
  independent, then by a theorem of Kummer theory (\cite{La}, Theorem
  8.1, p.~294) the square roots of these elements will generate a
  field extension of degree $2^s$. This extension is contained in
  $G$. So $\card(V)= 2^{\,\dim V}\leq r$.
\end{proof}


\end{document}